\documentclass[12pt]{amsart} 
 
\usepackage{amsthm}
\usepackage{amsmath}
\usepackage{amssymb}
\usepackage{pstricks, pst-node}

\usepackage{graphicx}
 
\newtheorem{theorem}{Theorem}
\newtheorem{ntheorem}{Theorem}
\newtheorem{definition}{Definition}

\newtheorem{nproposition}{Proposition}

\newtheorem{example}{Example}

\begin{document}

 \title{Threshold Graphs, Shifted Complexes, and Graphical
  Complexes}

 \author{Caroline J. Klivans}
 \address{ Depts. of Mathematics
  and Computer Science \\ The University of Chicago, Chicago, IL 60637}

\thanks{Author partially supported by NSF VIGRE grant DMS-0502215.}


\begin{abstract}

  We consider a variety of connections between threshold graphs,
  shifted complexes, and simplicial complexes naturally formed from a
  graph. These graphical complexes include the independent set,
  neighborhood, and dominance complexes. We present a number of
  structural results and relations among them including new
  characterizations of the class of threshold graphs.

\end{abstract}

\maketitle

 \section{Introduction}

 Threshold graphs are a well-studied class of graphs motivated from
 numerous directions.  They were first introduced by
 Chv\'atal and Hammer~\cite{Chvatal} as graphs for which there exists a linear
 threshold function separating independent from non-independent sets.
 Since then many equivalent conditions have been found for
 threshold graphs including constructive forms and forbidden
 configurations.  See for example~\cite{Mahadev} for nine different
 characterizations.

In generalizing to higher dimensions, it is then natural to consider
which characterizations remain equivalent. Golumbic first considered
such generalizations of threshold graphs to higher dimensions (or
hypergraphs)~\cite{Golumbic}.  He specifically highlighted three
analogs and asked if they were in fact the same. It turns out that
these three do not lead to the same class of
complexes~\cite{Reiterman}.  One of these analogs does give the class
known as shifted simplicial complexes.  We will primarily consider
threshold graphs from this perspective; that they are exactly the
one-dimensional shifted complexes.  See also~\cite{Klivans-Reiner}
which considers generalizations of threshold graphs based on degree
sequence properties and~\cite{Edelman} for a simple games/voting
theory perspective.

Shifted complexes are simplicial complexes whose faces form an order
ideal in the component-wise partial order. (See section~\ref{sec:def}
for precise definitions and examples of shifted complexes and
threshold graphs.)  Shifted complexes are named as such because of the
existence of shifting operations.  In general, a shifting operation
associates a shifted complex to any simplicial complex in a way which
preserves certain combinatorial properties but simplifies other
structure.  The original form of shifting, now known as combinatorial
shifting, was first introduced by Erd\"{o}s, Ko, and Rado~\cite{Erdos}
and Kleitman~\cite{Kleitman}.  More recently, Kalai~\cite{Kalai}
introduced algebraic shifting and spurred new interest in shifted
complexes.

 We study a variety of simplicial complexes naturally formed from a
 simple graph.  In many cases, these graphical complexes turn out to
 be shifted if and only if the graph is threshold.  We thus further
 motivate shiftedness as a natural generalization of threshold graphs.

 We first look at the independent set (or stable set) complex of a
 graph and show that it is shifted if and only if the graph is
 threshold.  Using this, we determine a constructible form for such
 complexes in terms of two simple operations.  Independent set
 complexes of graphs are also known as flag complexes. Combining this
 perspective and the construction, it is shown that pure shifted
 flag complexes are the same as pure shifted balanced complexes.

  Next we consider a generalized procedure to form the independent set
 complex of an arbitrary simplicial complex as in~\cite{Ehrenborg}.
 This construction again yields a shifted complex if and only if we
 start with a shifted complex.  Finally, we end with a result which
 shows that the dominance complex of a graph equals the neighborhood
 complex if and only if the graph is threshold.

\subsection{Definitions and Preliminaries}
\label{sec:def}

\begin{definition} A simplicial complex on $n$ vertices is {\em shifted} if there exists a
labeling of the vertices by one through $n$ such that for any face
$\{v_1, v_2, \ldots, v_k\}$, replacing any $v_i$ by a vertex with a
smaller label results in a collection which is also a face.
\end{definition}

An equivalent formulation of shifted complexes is in terms of order
ideals.  An order ideal $I$ of a poset $P$ is a subset of $P$ such
that if $x$ is in $I$ and $y$ is less than $x$ then $y$ is in $I$.
Let $P_s$ be the partial ordering on strings of increasing integers
given by $ { x} = (x_1 < x_2 < \cdots < x_k)$ is less than ${ y} =
(y_1 < y_2 < \cdots < y_k)$ if $x_i \leq y_i$ for all $i$ and ${ x}
\neq { y} $. Shifted complexes are exactly the order ideals of $P_s$.
We also allow comparisons of strings of various lengths by considering
the shorter string to have the necessary number of initial zeros
(slightly abusing that we are otherwise comparing strictly increasing
strings).  For example the string $24$ is taken to be less than the
string $1356$ by considering $24$ as $0024$.

\begin{figure}[h]
\centering
\includegraphics[height=2.7in, width=3.9in]{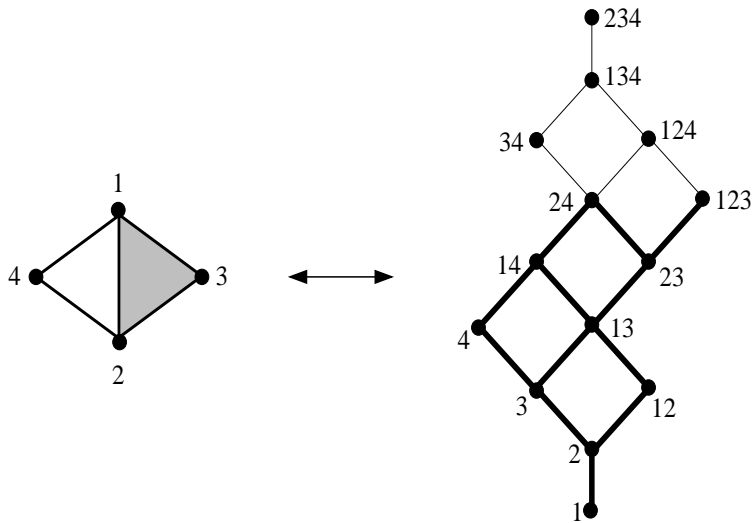}
\caption{An example of a shifted complex.}
\label{fig:1}
\end{figure}

\begin{example}
 \emph {A simplicial complex which includes the face
$\{24\}$ must also have the face $\{14\}$ in order to be shifted (see
Figure~\ref{fig:1}).}
\end{example}

One-dimensional shifted complexes are known to be the same as
threshold graphs~\cite{Kalai}. Threshold graphs are graphs that can be
given a vertex weighting which differentiates between independent and
non-independent sets.  An {\em independent set} of a graph is a collection
of vertices no two of which are connected.

\begin{definition}
\label{def:threshold}
A graph is threshold if for all $v \in V$ there exists weights $w(v)$,
and $t \in {\mathbb R}$ such that the following condition holds:
$w(U) \leq t$ if and only if $U$ is an independent set, where $w(U) =
\sum_{v \in U} w(v)$ (see Figure $2$).
\end{definition}

\begin{figure}[h]
\centering
\includegraphics[height = 1in]{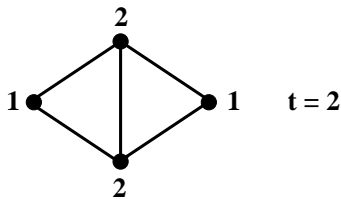}
\caption{A threshold graph with threshold 2.}
\label{fig:example2.4}
\end{figure}

One of the many characterizations of threshold graphs is constructive.
The construction is in terms of two basic operations; starring a
vertex and adding a disjoint vertex.  Starring a vertex $v$ onto a
graph $G = (V,E)$ forms the new graph:
$$G \, \, \textrm{star} \, \, v \, = \, (V \cup \{v\}, E \cup \{ \{x,v\} : x \in V\}).$$

Starring adds a new vertex adjacent to all previous vertices.

\begin{theorem}
[\cite{Mahadev} Theorem 1.2.4]
A graph is threshold if and only if it can be constructed from the one-vertex
graph by repeatedly adding a disjoint vertex or a starred vertex.
\end{theorem}

We want to extend the notion of starring a vertex to arbitrary
dimensions. Namely, we will say a vertex $v$ is starred in dimension
$d$ onto a complex $K$ by forming the complex: 
$$K \, \, \textrm{star}_d \, \, v \, = \, K \cup \{ v \cup f \, | \, f
\in K \, \, \textrm{and} \, \, |f| \leq d\}.$$ Note that this
operation is not the same as coning.  Coning corresponds to the
special case of starring a vertex in dimension one more than the
dimension of the complex. Coning will always increase the dimension of
a complex whereas starring does not necessarily increase the
dimension.  For example, let $K$ be the two-dimensional simplex
$\{123\}$.  $K$ star$_2$ $ 4$ is the two-dimensional boundary complex
of the three-dimensional simplex.  On the other hand, $K \, \,
\textrm{star}_3 \, \, 4$ is the three-dimensional simplex $\{1234\}$.

We will represent complexes generated by these two operations 
 as strings of $D$s (disjoint), $S$s (starring), and
vertical lines $|$ (for dimension increase).

\begin{example}  \rm \ Consider the string $DDSS|SSD|S$.  This
    represents the complex formed as follows: place two disjoint
    vertices, star two vertices in dimension $1$, star two vertices in
    dimension $2$, add a disjoint vertex, and star one vertex in
    dimension $3$.    Note that the string $DDSS|SS|DS$ would give the
    same complex.  For consistency, we will always place a vertical
    bar only immediately preceding an $S$.  Also note that it does not
    matter whether we begin a string with an $S$ or a $D$.  Again for
    consistency, we will always start a string with a $D$.  

    Given a complex represented by such a string, a shifted labeling
      can easily be obtained.  Suppose $K$ is a complex with $n$
      vertices, $k$ of which were added by starring.  Label the
      vertices corresponding to $S$ operations by $1$ through $k$ from
      right to left along the string.  Label the vertices represented
      by $D$ operations by $k+1$ through $n$ from left to right along
      the string.  The string above would give:

$$ \frac{D}{6} \frac{D}{7} \frac{S}{5}  \frac{S}{4} \frac{|}{} \frac{S}{3}
 \frac{S}{2}  \frac{D}{8} \frac{|}{} \frac{S}{1}.$$

\end{example}

While all complexes formed this way are shifted, not all shifted
complexes can be constructed by repeated application of these two
operations.  For example the complex of Figure~\ref{fig:1} does not have this form.

\section{Independence Complex of a Graph}

Recall that an {\em independent set} (or stable set) of a graph is a
collection of vertices no two of which are connected by an edge.  Let
$I(G)$ denote the independence complex of a graph $G$.  This complex
is formed by taking the collection of independent sets of $G$.
Clearly removing a node from an independent set results in an
independent set so this collection is a simplicial complex.

\begin{ntheorem}
\label{Thm:I(G)}
 $I(G)$ is shifted if and only if $G$ is a threshold graph. 
\begin{proof}
Let $G$ be a threshold graph.  Then we know $G$ is shifted.  Let $l$
be a shifted labeling of the vertices of $G$.  Consider any face
$F = \{v_1, v_2 \ldots v_k \}$ of $I(G)$ and a vertex $w$ such that $l(w) >
l(v_i)$ for some $i$.  We will show that $ F' = \{v_1, v_2, \ldots,
\hat{v_i},w, \ldots, v_k \}$ is a face of $I(G)$.  If not, then $w$ 
must be connected to some $v_j \, (j \neq i)$ in $G$. Because $w$ has
a larger label than $v_i$ and  $\{wv_j\} \in E(G)$, $\{v_iv_j\}$ must be an edge of $G$ in order
for $G$ to be shifted.  But this contradicts $F$ being a face of $I(G)$.  Hence $I(G)$ is shifted under the reverse ordering of $l$.

Now let $I(G)$ be shifted and $l$ a shifted labeling.  Consider any edge, 
$\{v_1v_2\}$ of $G$ and a vertex $w$ such that $l(w) > l(v_2)$.  We will
 show that $\{v_1w\}$ is an edge of $G$ and hence $G$ is shifted again under the reverse ordering of $l$.  If not, then $\{v_1w\}$ is an 
independent set of $G$ and hence a face of $I(G)$.  $I(G)$ is shifted and 
$v_2$ had a smaller label than $w$ which means $\{v_1v_2\}$ must be a face of $I(G)$
and not an edge of $G$, again a contradiction.
  
\end{proof}
\end{ntheorem}

\subsection{Flag complexes}

Independent set complexes of graphs are also known as flag complexes.
A {\em flag} complex is defined as a simplicial complex such that
every minimal non-face has exactly two elements~\cite{Stanley}. By the previous
result, all shifted flag complexes are formed from threshold graphs.
Using both perspectives allows us to further determine the form of
these complexes.

\begin{ntheorem}
Shifted flag complexes are the complexes formed by the operations $D$ and $S$ with exactly one $S$ in each dimension.
\end{ntheorem}
\begin{proof}

Every shifted flag complex arises as the independence complex of a
threshold graph.  Every threshold graph can be represented as a
string of $D$s and $S$s.  Consider mapping this string under the following
rules: $D \rightarrow |S$ and $S \rightarrow D$.  Namely, switch every
$S$ to a $D$ and switch every $D$ to an $S$ and also increase the
dimension with every such switch.

\begin{example}  $DDSDSDSSD \rightarrow S|SD|SD|SDD|S$
\end{example}

First, we want to determine the independent sets of a threshold graph from its
string of $D$s and $S$s.  The maximal independent sets are the set of all $D$s and all collections
which consist of a single $S$ and all $D$s that come after it.  

Next, given the image of the string, we want to determine its facets.  They are
the set of all $S$s and all collections which consist of a $D$ and all $S$s that 
come after it.  In particular they are exactly the independent sets of $G$.

This procedure is invertible showing that all strings of $D$s and $S$s with exactly
one $S$ in each dimension are flag complexes.  

\end{proof}

\subsection{Balanced complexes}
\label{sec:balanced}

A $d$-dimensional simplicial complex is {\em balanced} if its
 vertices can be colored with $d+1$ colors such that within any face
 all vertices have different colors.

\begin{nproposition}
All shifted flag complexes are balanced.
\end{nproposition}

\begin{proof}
  Let $K$ be a $d$-dimensional shifted flag complex. We give an
  explicit balanced labeling.  $K$ can be represented as a string of
  $D$s and $S$s with exactly one star operation per dimension.  Label
  the vertices with $d+1$ colors as shown below:

\[ \underbrace{DD \dots D}_{1} \underbrace{SD \ldots
  D}_{2}| \ldots |\underbrace{SD \ldots D}_{d}|\underbrace{SD \ldots
  D}_{d+1}. \]

Every face of $K$ consists of one initially placed disjoint vertex and
a set of starred vertices which come after it, all of which have been
given a different color.

\end{proof}

The converse of the proposition above is false: not all shifted balanced complexes are flag
complexes.  A simple example is the complex on 4 vertices with maximal
faces $\{123, 14, 24\}$ (see Figure~\ref{fig:1}).  Notice that this complex is not pure. 

A pure shifted flag complex has a very simple form:
$$DD \ldots DS|S|S \ldots |S|S.$$  This yields a ``pencil of facets''.
Namely, a $d$-dimensional pure shifted flag complex on $n$ vertices
consists of $n-d$ facets all sharing a common $d-1$ face.

\begin{ntheorem}
A pure shifted complex is balanced if and only if it is a flag complex.
\end{ntheorem}

\begin{proof}
  We already know flag implies balanced.  We will show any pure
  shifted balanced complex is also a ``pencil of facets''.  Let $K$ be
  a $d$-dimensional pure, shifted, and balanced complex with a shifted
  labeling of its vertices.  Shiftedness implies that $\{1,2, \ldots,
  d+1\}$ $\in$ $K$.  Let $\{x_1,x_2, \ldots, x_{d+1}\}$ ($x_1 < x_2
  \cdots <x_{d+1}$) be some other facet.  Suppose $x_{d} > d$.  Then
  $x_{d+1}$ must be greater than $d+1$ and vertex ${d+2}$ must be
  adjacent to $d+1$.  But then shiftedness implies that the complete
  graph on $d+2$ vertices is in the $1$-skeleton of $K$ which
  contradicts $K$ being balanced.  Hence $x_d$ must equal $d$ and all
  facets have the form $\{1,2, \ldots, d, x\} \, \, d < x \leq n$.  Thus
  $K$ has the same form as a pure shifted flag complex. 

\end{proof}

\subsection{Shifting}

Recall that in general, a shifting operation associates a shifted
complex to any simplicial complex in a way which preserves certain
combinatorial properties but simplifies other structure.  In
particular, both combinatorial and algebraic shifting preserve the
$f$-vector.  

A conjecture due to Kalai asks if any $f$-vector of a flag complex can
also be realized as the $f$-vector of a balanced
complex~\cite{Stanley}.  We note here for completeness the
relationship between shifting and the properties of flag and
balanced.  The two main variants of algebraic shifting unfortunately
do not preserve flag or balanced complexes. For definitions and much
more on these shifting operations see~\cite{Kalai}.  Consider the
complete bipartite graph $K_{3,3}$. It is easy to check that this is
both a flag and balanced complex.  Symmetric shifting yields the
complex generated by top face $\{26\}$ and exterior shifting yields
the complex generated by top faces $\{25\}$ and $\{34\}$.  In both
graphs, the collection $\{123\}$ is a minimal non-face showing it is
not a flag complex and not balanced.

Moreover, no shifting operation which preserves the $f$-vector
could preserve these properties.  The graph $K_{3,3}$ has 6 vertices,
9 edges, and no faces of dimension 2 or greater.  But any order ideal
in the shifted partial order on 6 vertices with 9 one-dimensional
faces will include the edges $\{12\}$, $\{13\}$, and $\{23\}$.  Hence
the graph will not be balanced and since we can not add any
two-dimensional faces, this will generate a minimal non-face with
three elements.

\section{Generalized Independence Complex}

In~\cite{Ehrenborg}, forming the independence complex of a graph is
 generalized to arbitrary simplicial complexes.  For a simplicial complex $K$, define $I(K)$
by declaring the facets of $K$ to be the minimal non-faces of $I(K)$.
(The independent set complex in~\cite{Ehrenborg} is defined in greater
generality, allowing for $K$ to be a set system which is not
necessarily a simplicial complex.)

We start by considering the independent set complex of shifted
simplicial complexes.  The general statement that $K$ is shifted if
and only if $I(K)$ is shifted is false in both directions.  It is not
hard to construct counter-examples using non-pure complexes.  For
example, let $K$ be the simplicial complex on $5$ vertices with
maximal faces $\{123, 14, 24, 15\}$.  $K$ is shifted but $I(K)$ which
has maximal faces $\{235,345,12,13\}$ is not.  The induced subcomplex
on vertices $\{1,2,4,5\}$ is a path of length three which is an
obstruction to shiftedness in dimension one.  We can continue to apply
the procedure to disprove the other direction.  $I(I(\Delta))$ is
generated by $\{245,234,145,35\}$ which is also not shifted.
$I(I(I(\Delta)))$ generated by $\{123,124,125,134,45\}$ is however
shifted (mapping $3 \leftrightarrow 4$ gives a shifted labeling).

Restricting to the pure case is actually a more natural generalization
of the independence complex of a graph.  The generalized procedure
only restricts to the same procedure on graphs if the graph is
connected (i.e. pure).  For example, if we have a graph with disjoint
vertices, under the generalized procedure they would be minimal
non-faces of $I(K)$.  On the other hand, a disjoint vertex is in all
maximal faces of the independence complex of the graph.  In the pure
case, we come to the following result:

\begin{ntheorem}
For $K$ pure, 
$K$ is shifted if and only if $I(K)$ is shifted.
\end{ntheorem}

\begin{proof}

 Suppose $K$ is shifted but $I(K)$ is not shifted.  Then there exists
$x,y,f_1,f_2 \in I(K)$ such that $xf_1, yf_2 \in I(K)$ and $yf_1, xf_2
\notin I(K)$, where $x$ and $y$ are vertices and $f_1$ and $f_2$ are
faces, see~\cite{Klivans}.  Since $yf_1$ and $xf_2$ are not in $I(K)$,
they must be facets or contain facets of $K$.  First we note that the
facets involved here must not be strictly contained in $f_1$ and $f_2$,
or $xf_2$ and $yf_1$ could not be in $I(K)$.

Suppose $yf_1$ and $xf_2$ are facets of $K$.  Let $l$ be a
shifted labeling for $K$ and without loss of generality, let
$l(x) < l(y)$.  Since $K$ is shifted, we have that $xf_1 \in
K$. But, $|xf_1|$ = $|yf_1|$ which implies $xf_1$ is a facet of
$K$ and can not be in $I(K)$ - a contradiction.

Suppose at least one of $yf_1$ and $xf_2$ is not a facet of $K$.
They still must contain a facet.  Let $g_1 \subseteq f_1$, $g_2
\subseteq f_2$, and $xg_2$, $yg_1$ be facets of $K$.  They will
not be in $I(K)$, but $yg_2 \subseteq yf_2 \in I(K)$ and
$xg_1 \subseteq xf_1 \in I(K)$ so we are back in the first case.

Now suppose $I(K)$ is shifted but $K$ is not shifted.
Then there exists $x,y,f_1,f_2$ such that $xf_1, yf_2 \in K$ and
$yf_1, xf_2 \notin K$. Because $K$ is pure, we may take $xf_1$ and
$yf_2$ to be maximal faces; in particular this gives that $|xf_1|$ =
$|yf_2|$.  Now since $xf_1$ and $yf_2$ are facets of $K$, they are not
in $I(K)$.  Next consider $xf_2$ and $yf_1$.  For these faces not to
be in $I(K)$, they must contain facets.  However, $|xf_2|$ = $|yf_2|$
= $|xf_1|$ = $|yf_1|$ so if they contained a facet it would be of
smaller size, and this can not be because $K$ is pure.  Hence $xf_2$
and $yf_1$ are in $I(K)$, which contradicts $I(K)$ being shifted.

\end{proof}

\subsection{Neighborhood and Dominance}

A {\em dominating} set of a graph is a set of vertices $D$ such that
all vertices are either in $D$ or adjacent to a vertex in $D$.  The
{\em dominance complex} $D(G)$ is the collection of {\bf complements}
to dominating sets~\cite{Ehrenborg}. Note that this is because dominating sets are closed under
superset as opposed to subsets.

In~\cite{Ehrenborg} the dominance complex is studied for specific
graphs and it is observed that $D(G)$ is the independent set complex of
the collection of closed neighborhoods $N[v]$ of $G$. (The closed
neighborhood of a vertex $v$ is the usual neighborhood $N(v)$ union
$v$ itself). If we define the {\em closed neighborhood complex} $N[G]$ to be the
simplicial complex with facets equal to the minimum (under inclusion)
sets of the collection of closed neighborhoods of $G$, then $I(N[G])$
as we have defined $I(K)$ matches that of~\cite{Ehrenborg}.

By the previous result, we would hope to show that $G$ is threshold
if and only if $N[G]$ is shifted, and hence $D(G) = I(N[G])$ is shifted if and only if $G$
is threshold.  This is unfortunately not the case.

We do however offer a curious relationship between these and the usual
neighborhood complex of Lovasz~\cite{Lovasz}.  Let $N(G)$ be the
collection of sets of vertices which share a common neighbor.

\begin{ntheorem}
$N(G) = D(G)$ (and therefore $I(N[G]$)) if and only if $G$ is threshold
\end{ntheorem}

\begin{proof}
  First we claim that $N(G) \subseteq D(G)$ for any graph.  Suppose
  $\{x_1, x_2, \ldots, x_k\}$ $\in N(G)$.  Let $v$ be their common
  neighbor.  Now $v \in V \, \backslash \,  \{x_1,x_2, \ldots, x_k\}$ so all
  vertices are either in the complement or adjacent to a vertex in the
  complement, hence it dominates.  

  Next we show that $D(G) \subseteq N(G)$ if $G$ is threshold.  Let
  $G$ be threshold with a shifted labeling $l$ and let $\{x_1, x_2,
  \ldots, x_k\} \in D(G)$.  We need to show that the $ \{x_i\}$ have a
  common neighbor.  Without loss of generality let $x_k$ have the
  largest label among the $x_i$.  Because $V \, \backslash \,  \{x_1, x_2, \dots,
  x_k\}$ dominates, every $x_i$ is adjacent to some vertex in $V \, \backslash \,   \{x_1, x_2, \dots, x_k\}$.  Let $x_k$ be adjacent to some vertex
  $v$.  Because $x_k$ has the largest label and $G$ is threshold, $v$
  is a common neighbor to all the $x$s.

  Finally, we show that if $N(G) = D(G)$ then $G$ is threshold.  We
  will do this by showing that $G$ is constructed by repeatedly
  starring or adding a disjoint vertex.  Let $G$ be such that $N(G) =
  D(G)$.

  First, at most one connected component of $G$ has an edge.  Suppose
  more than one connected component had an edge.  The complement to
  any minimal dominating set must contain vertices from different
  connected components.  Hence the complement can not have a common
  neighbor.  Note that a component which is a single vertex is fine,
  this vertex will be in all dominating sets.  

  Next consider the connected component with at least one edge, if
  there is no such component then the graph is a collection of
  disjoint vertices which is shifted.  Otherwise, we claim it has a
  star vertex.  Suppose not, then any minimal dominating set has size
  at least two.  Let $D =\{x_1, x_2, \ldots, x_k\}$ $(k \geq 2)$ be a
  dominant set of minimal size.  Then $V \, \backslash \, D$ is a
  maximal element of $N(G)$ and hence the neighborhood, $N(v)$, of
  some vertex $v$.  Note that $v \in D$ or else both $v$ and $N(v)$
  are in the complement and $D$ could not be dominating.  Without loss
  of generality let $v = x_1$.  Consider another vertex $y$ such that
  $x_1y$ and $x_2y$ are edges of $G$ where $x_2$ has also been taken
  without loss of generality.  Such a $y$ must exist because we are
  working in a connected component and the $x_i$s can not be adjacent
  to each other because $D$ is minimal.  $y$ can not be in an edge
  with any other $x_i$ or else $D \, \backslash \, \{x_2, x_i\} \cup
  \{y\}$ would be a smaller dominating set.  Therefore $D' = D \,
  \backslash \, \{x_2\} \cup \{y\}$ is another minimal dominating set.
  Hence the complement of $D'$ must be the neighborhood of some
  vertex, say $w$.  Now, $w$ can not equal any $x_i$ or else
  $\{x_2x_i\}$ is an edge which contradicts the minimality of $D$.
  And, $w \neq y$ or else $D' \, \backslash \, x_1$ would be
  dominating which contradicts the minimality of $D'$.  Hence we've
  reached a contradiction since $w \in D'$ must hold.  Therefore the
  connected component of $G$ with at least one edge has a star vertex.

  To finish the proof, we only need to show that if $N(G) = D(G)$ then
  $N(G \, \backslash \, v) = D(G \, \backslash \, v)$ for $v$ a star
  vertex.  (Removing any disjoint vertices does not affect either
  complex).  Clearly, $N(G \, \backslash \, v)$ $= \{f \, \backslash
  \, v \, | \, f \in N(G), v \in f\}$.  Note that the only maximal
  face of $N(G)$ which does not contain $v$ is $V \, \backslash \, v$.

  Similarly, moving from $D(G)$ to $D(G \, \backslash \,  v)$ we lose the one
  facet of $D(G)$ corresponding to all vertices except $v$.  Any
  minimal dominating set for $G$ (other than the set $\{v\}$) is
  dominating for $G \, \backslash \,  v$ as well.  Therefore $D(G \, \backslash \,  v)$ also
  equals $ \{f \, \backslash \,  v \, | \, f \in D(G), v \in f\}$.

  Because threshold graphs are exactly those graphs which can be
  constructed by repeatedly adding a disjoint vertex and a star vertex,
  $G$ is threshold.

\end{proof}

\end{document}